\documentclass[11pt,a4paper]{article}

\usepackage{array}
\newcolumntype{"}{@{\hskip\tabcolsep\vrule width 1.5pt\hskip\tabcolsep}}

\usepackage{graphics}
\usepackage{epsfig}
\usepackage{subcaption}
\usepackage{enumerate}
\usepackage{amssymb}
\usepackage{amsthm}
\usepackage{mathrsfs}
\usepackage{array}
\usepackage{hhline}
\usepackage{longtable}
\usepackage{amsmath}
\usepackage{float}
\usepackage{picinpar}
\usepackage{cite}
\usepackage{xcolor}
\usepackage{multirow}
\usepackage[mathlines]{lineno}
\modulolinenumbers[2]
\newcommand*\patchAmsMathEnvironmentForLineno[1]{%
\expandafter\let\csname old#1\expandafter\endcsname\csname #1\endcsname  \expandafter\let\csname oldend#1\expandafter\endcsname\csname end#1\endcsname  \renewenvironment{#1}%
{\linenomath\csname old#1\endcsname}%
{\csname oldend#1\endcsname\endlinenomath}}%
\newcommand*\patchBothAmsMathEnvironmentsForLineno[1]{%
\patchAmsMathEnvironmentForLineno{#1}%
\patchAmsMathEnvironmentForLineno{#1*}}%
\AtBeginDocument{%
\patchBothAmsMathEnvironmentsForLineno{equation}%
\patchBothAmsMathEnvironmentsForLineno{align}%
\patchBothAmsMathEnvironmentsForLineno{flalign}%
\patchBothAmsMathEnvironmentsForLineno{alignat}%
\patchBothAmsMathEnvironmentsForLineno{gather}%
\patchBothAmsMathEnvironmentsForLineno{multline}%
}

\definecolor{vividviolet}{rgb}{0.62, 0.0, 1.0}

\def\rsq{\hspace*{\fill}$\blacksquare$\medskip}

\newtheorem{theorem}{Theorem}[section]
\newtheorem{lemma}[theorem]{Lemma}
\newtheorem{corollary}[theorem]{Corollary}
\newtheorem{problem}{Problem}[section]

\newtheorem{conjecture}{Conjecture}[section]

\numberwithin{equation}{section}
\newtheorem{observation}[theorem]{Observation}

\newtheoremstyle{example}
  {10pt}          
  {10pt}  
  {\rm}  
  {}
  {\bf}  
  {: }    
  { }    
  {}     
\theoremstyle{example}
\newtheorem{example}{Example}[section]

\newtheorem{defi}{Definition}[section]

\usepackage{graphicx}
\graphicspath{%
    {converted_graphics/}
    {/}
}

\def\ms{\medskip}

\def\nt{\noindent}

\begin{document}

\begin{center}
{\mathversion{bold}\Large \bf Sudoku Number of Graphs}

\bigskip
{\large  Gee-Choon Lau$^{a,}$\footnote{Corresponding author.}, J. Maria Jeyaseeli$^{b}$, Wai Chee Shiu$^{c}$, S. Arumugam{$^{b}$}}\\

\medskip

\emph{{$^a$}Faculty of Computer \& Mathematical Sciences,}\\
\emph{Universiti Teknologi MARA (Johor Branch, Segamat Campus),}\\
\emph{85000, Malaysia.}\\
\emph{geeclau@yahoo.com}\\

\ms

\emph{{$^b$}National Centre for Advanced Research in Discrete Mathematics,\\ Kalasalingam Academy of Research and Education,\\
Anand Nagar, Krishnankoil 626 126, Tamil Nadu, India}\\
\emph{mariyajeyaseeli@gmail.com,s.arumugam.klu@gmail.com}\\

\medskip
\emph{{$^c$}Department of Mathematics,}\\
\emph{The Chinese University of Hong Kong,}\\
\emph{Shatin, Hong Kong, China.}\\
\emph{wcshiu@associate.hkbu.edu.hk}\\
\end{center}

\begin{abstract}
We introduce a new concept in graph coloring motivated by the popular Sudoku puzzle. Let $G=(V,E)$ be a graph of order $n$ with chromatic number $\chi(G)=k$ and let $S\subseteq V.$ Let $\mathscr C_0$ be a $k$-coloring of the induced subgraph $G[S].$ The coloring $\mathscr C_0$ is called an extendable coloring if $\mathscr C_0$ can be extended to a $k$-coloring of $G.$ We say that $\mathscr C_0$ is a Sudoku coloring of $G$ if $\mathscr C_0$ can be uniquely extended to a $k$-coloring of $G.$ The smallest order of such an induced subgraph $G[S]$ of $G$ which admits a Sudoku coloring is called the Sudoku number of $G$ and is denoted by $sn(G).$ In this paper we initiate a study of this parameter. We first show that this parameter is related to list coloring of graphs. In Section 2, basic properties of Sudoku coloring that are related to color dominating vertices, chromatic numbers and degree of vertices,  are given. Particularly, we obtained necessary conditions for $\mathscr C_0$ being uniquely extendable, and for $\mathscr C_0$ being a Sudoku coloring. In Section 3, we determined the Sudoku number of various familes of graphs. Particularly, we showed that a connected graph $G$ has $sn(G)=1$ if and only if $G$ is bipartite. Consequently, every tree $T$ has $sn(T)=1$. Moreover, a graph $G$ with small chromatic number may have arbitrarily large Sudoku number. Extendable coloring and Sudoku coloring are nice tools  for providing a $k$-coloring of $G$.

\noindent {\bf Keywords:} chromatic number, extendable coloring, Sudoku coloring.

\noindent {\bf 2010 AMS Subject Classifications:} 05C78; 05C69.
\end{abstract}

\section{Introduction}

By a graph $G=(V,E)$ we mean a finite undirected connected simple graph of order $n=|V|$ and size $m=|E|.$ For notations and concepts not defined here we refer to the book \cite{Bondy}.

\ms\nt  We introduce a new concept in graph coloring which is motivated by the popular Sudoku puzzle. Let $V$ denote the 81 cells consisting of 9 rows, 9 columns and nine $3\times 3$ subsquares. Let $G$ be the graph with vertex set $V$ in which each of the rows, columns and the nine $3\times 3$ subsquares are complete graphs. Clearly $\chi(G)=9$ and the solution of a Sudoku puzzle gives a proper vertex coloring of $G$ with nine colors. A Sudoku puzzle corresponds to a proper vertex coloring $\mathscr C$ of an induced subgraph $H$ of $G$ using at most 9 colors with the property that $\mathscr C$ can be uniquely extended to a 9-coloring of $G$. This motivates the following definition.

\ms
\begin{defi} Let $G$ be a connected graph of order $n$ with chromatic number $\chi(G)=k\ge 2$. Let $S \subset V$ and let $G[S]$ be the induced subgraph of $G.$  Let $\mathscr C_0$ be a proper $k$-coloring of $G[S]$.  We also call $\mathscr C_0$ a partial coloring of $G$. The coloring $\mathscr C_0$ is called an {\it extendable} coloring of $G[S]$ if $\mathscr C_0$ can be extended to a $k$-coloring for $G$. We say $\mathscr C_0$ is a {\it Sudoku coloring} of $G$ if $\mathscr C_0$ can be uniquely extended to a proper $k$-coloring of $G$. The smallest order of such an induced subgraph $G[S]$ of $G$ which admits a Sudoku coloring is called the {\it Sudoku  number} of $G$ and is denoted by $sn(G)$.
\end{defi}

\nt It is known that Arto Inkala provided a Sudoku puzzle with 21 initial numbers. Thus, the Sudoku number of the traditional Sudoku puzzle is at most 21 (see~\cite{Sudoku}).

\ms\nt We need the following definitions.

\begin{defi}[\hspace*{-1mm}\cite{Bondy}]  Let $G=(V,E)$ be a graph with chromatic number $k$. The graph $G$ is uniquely colorable if the partition $\{V_1,V_2,\dots,V_k\}$ of $V$ induced by any $k$-coloring of $G$ is unique up to a permutation.
\end{defi}

\nt The concept of list coloring was independently studied by Erd\H{o}s {\it et al.} \cite{Erdos} and Vizing \cite{Vizing}.

\begin{defi} Given a graph $G$ and a set $L(v)$ of colors for each vertex $v,$ which is called a {\it list} or {\it color-list} of $v$, a list coloring is a coloring of $G$ such that the vertex $v$ is colored by a color in the list $L(v).$ If $G$ admits a list coloring for a given list $L,$ then $G$ is called {\it list colorable} or {\it $L$-colorable}.
\end{defi}

\ms\nt Scheduling problem is one of the classical applications of graph colorings. If the graph associated with an instance of the problem has chromatic number $k$, then $k$ is the minimum number of timeslots required for a scheduling and each color class represents a time slot.  Thus for implementation of a scheduling, the knowledge of the chromatic number alone is not sufficient and we need a $k$-coloring of $G$.  Extendable coloring and Sudoku coloring are nice tools  for providing a $k$-coloring of $G$.

\ms\nt In this paper we formulate Sudoku coloring as a list coloring problem and determine the Sudoku number of several classes of graphs.

\section{Extension: A List Coloring Problem}

\nt Let $\mathscr C_0$ be an extendable coloring of $G$. Then for each uncolored vertex $v$ there are some possible colors that can be used to color the vertex $v$. The set of all possible colors for $v$ is a color-list of $v$. The problem of extending $\mathscr C_0$ to a $k$-coloring of $G$ is equivalent to a list coloring problem.

\ms \nt In this paper, we use $\{i\;|\; 1\le i\le k\}$ as the color set for a graph $G$ if $\chi(G)=k$.

\begin{figure}[H]
\centering
\captionsetup[subfigure]{labelformat=simple}
\renewcommand{\thesubfigure}{(\alph{subfigure})}
\begin{subfigure}{0.3\textwidth}
\centering
{\epsfig{file=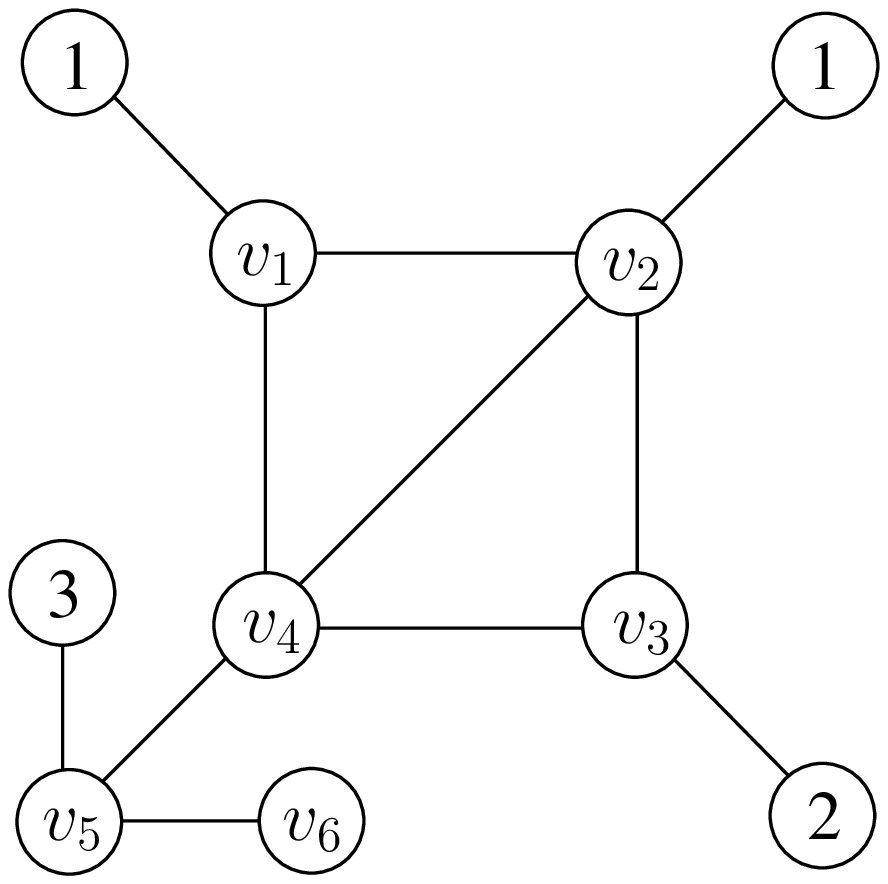, width=3cm}}
\caption{}\label{fig-ex0-1}
\end{subfigure}
\begin{subfigure}{0.3\textwidth}
\centering
\epsfig{file=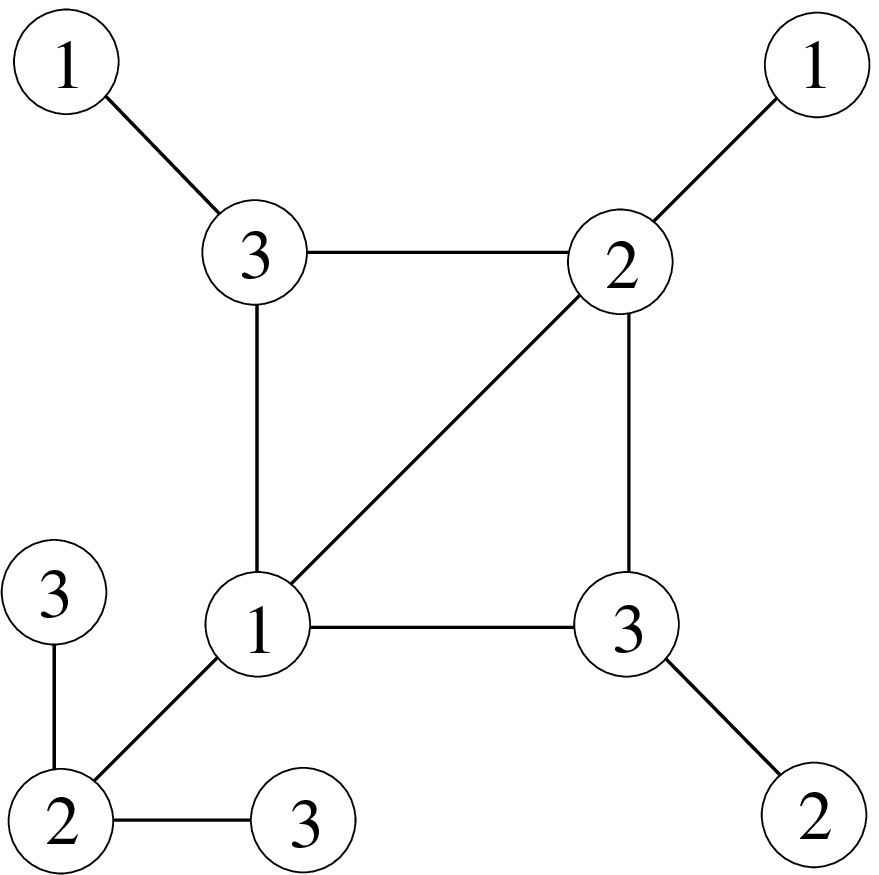, width=3cm}
\caption{}\label{fig-ex0-2}
\end{subfigure}
\begin{subfigure}{0.3\textwidth}
\centering
\epsfig{file=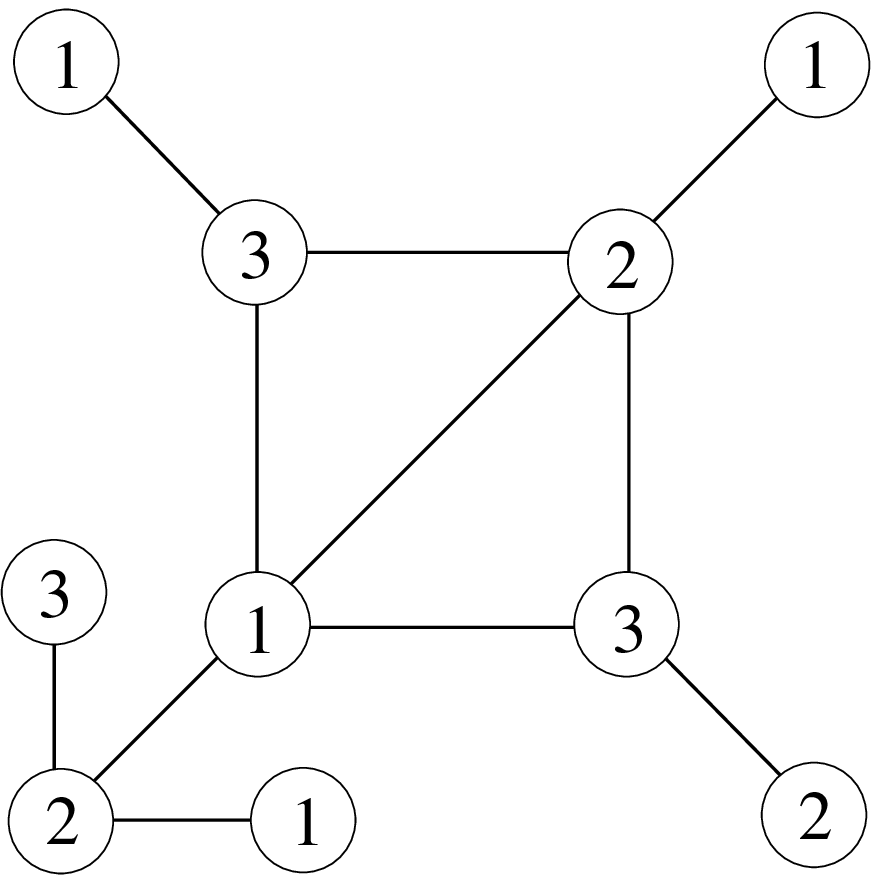, width=3cm}
\caption{}\label{fig-ex0-3}
\end{subfigure}
\caption{An extendable coloring with at least two extensions.}
\end{figure}

\begin{example} The partial coloring $\mathscr C_0$ for the graph given in Fig.~\ref{fig-ex0-1} is an extendable coloring. Also two extensions of $\mathscr C_0$ to a 3-coloring of $G$ are given in Fig.~\ref{fig-ex0-2} and Fig.~\ref{fig-ex0-3}. Thus $\mathscr C_0$ is an extendable coloring, but not a Sudoku coloring.
\end{example}

\rsq

\begin{example}  Consider the partial coloring $\mathscr C_0$ for the graph given in Fig.~\ref{f2}.
\begin{figure}[H]
\centering
\epsfig{file=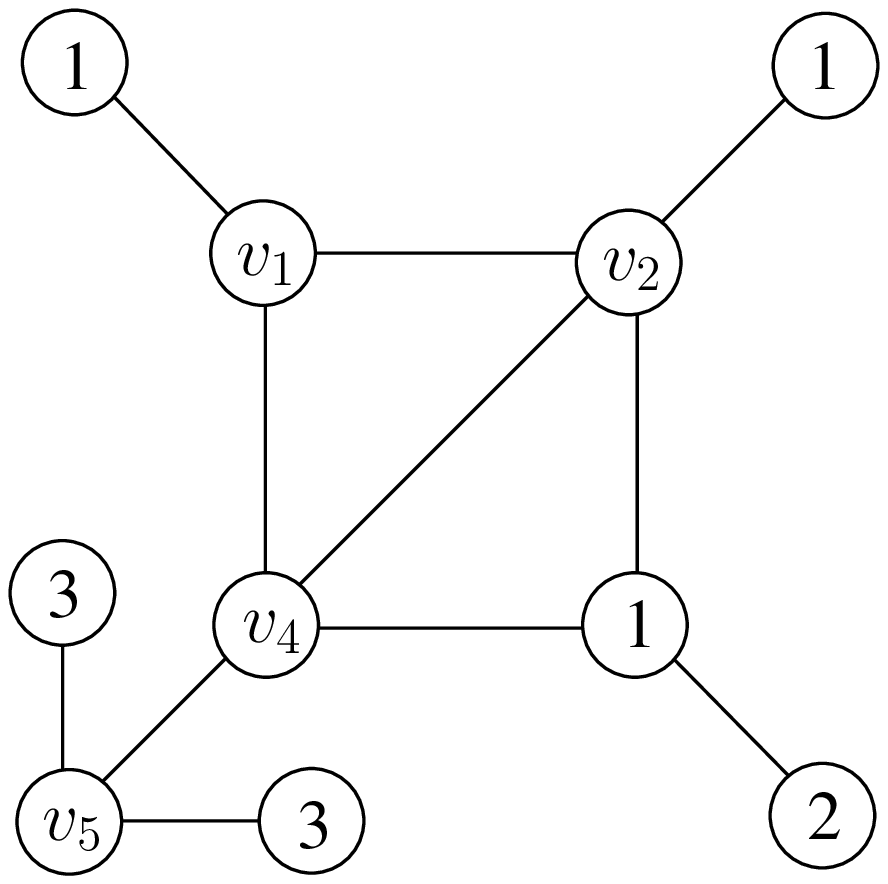, width=3cm}
\caption{A coloring which is not extendable.}\label{f2}
\end{figure}

\nt Clearly $L(v_4)=L(v_1)=L(v_2)=\{2,3\}$ and $G[\{v_1,v_2,v_4\}]=K_3.$ Hence $G$ is not $L$-colorable and so $\mathscr C_0$ is not an extendable coloring. \rsq
\end{example}

\begin{example}\label{ex-ex1} Consider the partial coloring $\mathscr C_0$, which involves 2 colors, for the graph given in Fig.~\ref{fig-ex1}.

\begin{figure}[H]
\centering
\captionsetup[subfigure]{labelformat=simple}
\renewcommand{\thesubfigure}{(\alph{subfigure})}
\begin{subfigure}{0.2\textwidth}
\centering
{\epsfig{file=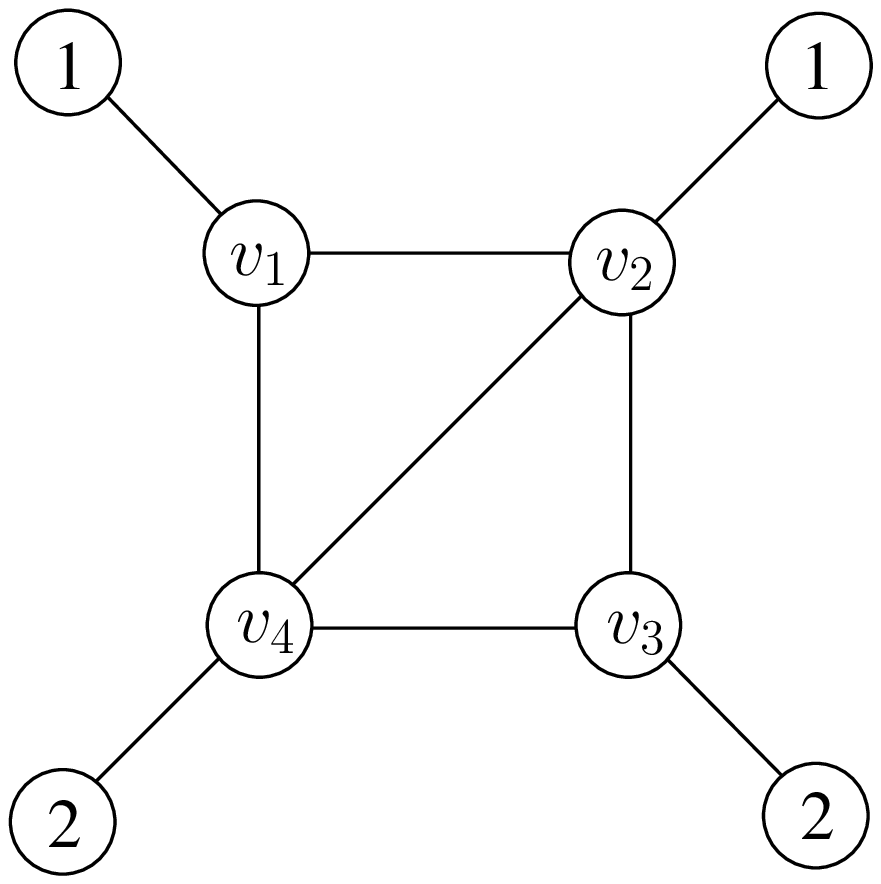, width=3.6cm}}
\caption{Coloring $\mathscr C_0$}\label{fig-ex1}
\end{subfigure}
\begin{subfigure}{0.37\textwidth}
\centering
\raisebox{-1.8cm}[1.5cm][2.1cm]{\epsfig{file=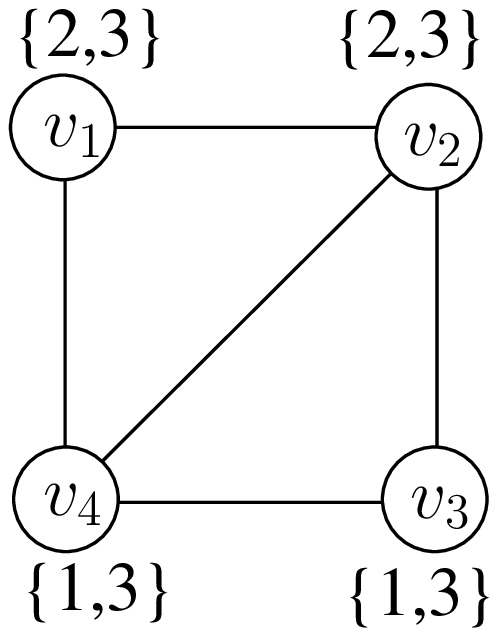, width=2cm}}
\caption{Color-lists $L$ for the uncolored vertices}\label{fig-ex1list}
\end{subfigure}
\begin{subfigure}{0.2\textwidth}
\centering
\epsfig{file=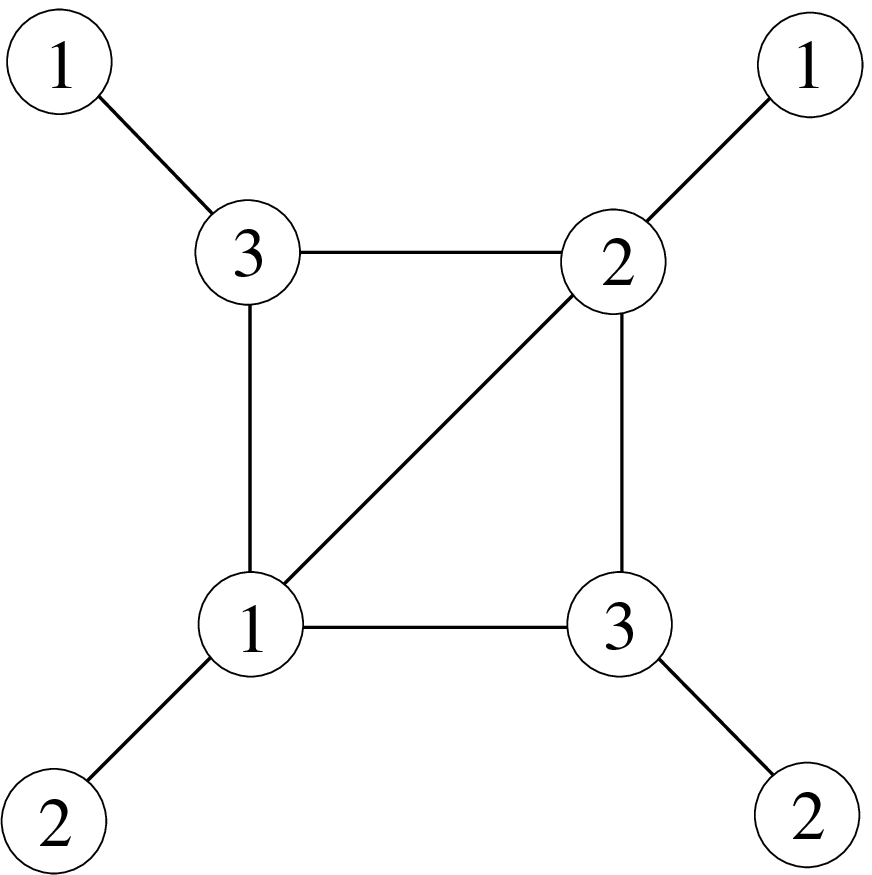, width=3.6cm}
\caption{Coloring $\mathscr C$}\label{fig-ex2}
\end{subfigure}
\caption{A Sudoku coloring and its unique extension.}\label{f3}
\end{figure}

\nt Clearly $L(v_2)=L(v_1)=\{2,3\}$ and $L(v_3)=L(v_4)=\{1,3\}$.  So we have Fig.~\ref{fig-ex1list}.  Let $\mathscr C$ be an extension of $\mathscr C_0$ to  a 3-coloring of $G.$ If  $v_2$ is assigned color 3, then the color list of the two adjacent vertices $v_3$ and $v_4$ becomes $\{1\}.$ Hence the color for $v_2$ is 2. Now the colors to be assigned to the remaining vertices are uniquely determined and hence the coloring $\mathscr C$ given in Fig.~\ref{fig-ex2} is the unique extension of $\mathscr C_0$. Thus $\mathscr C_0$ is a Sudoku coloring of $G$. \rsq
\end{example}

\begin{example}\label{ex-exI} The partial coloring $\mathscr C_0$, which involves 3 colors, for the graph given in Fig.~\ref{fig-exI} is a Sudoku coloring.

\begin{figure}[H]
\centering
\captionsetup[subfigure]{labelformat=simple}
\renewcommand{\thesubfigure}{(\alph{subfigure})}
\begin{subfigure}{0.2\textwidth}
\centering
{\epsfig{file=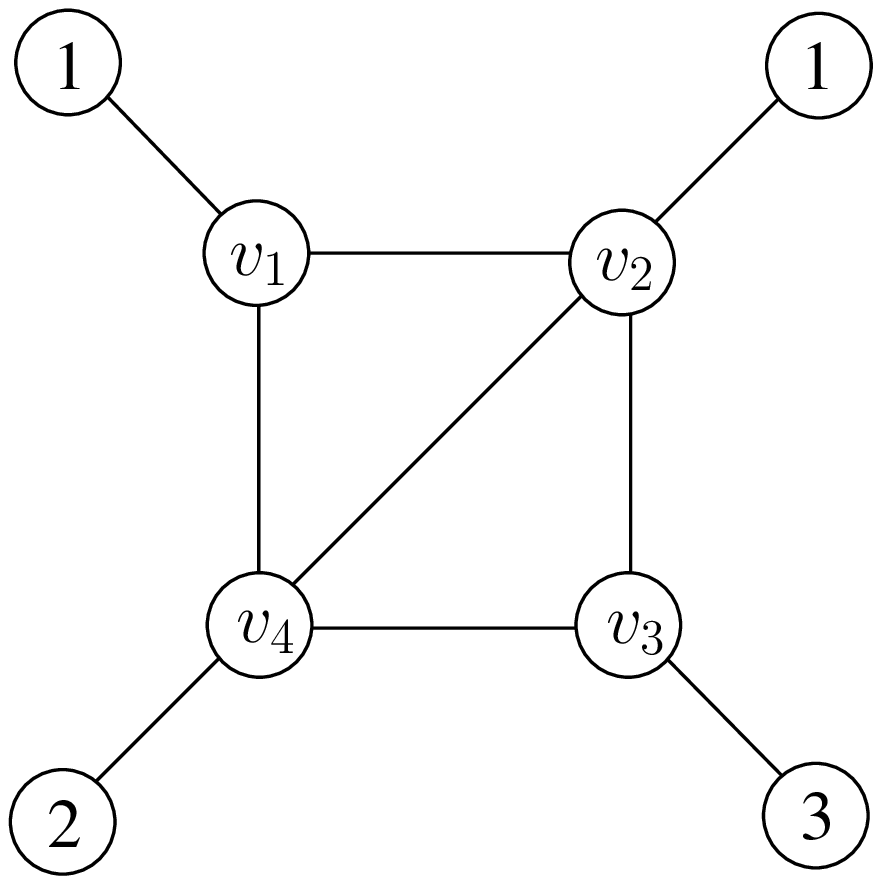, width=3.6cm}}
\caption{Coloring $\mathscr C_0$}\label{fig-exI}
\end{subfigure}
\begin{subfigure}{0.37\textwidth}
\centering
\raisebox{-1.8cm}[1.5cm][2.1cm]{\epsfig{file=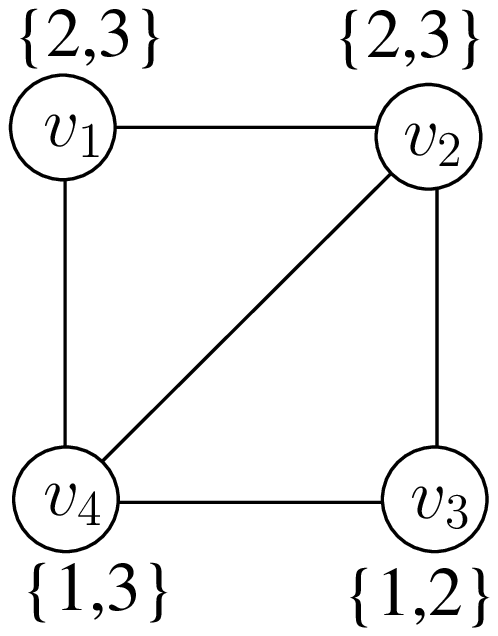, width=2cm}}
\caption{Color-lists $L$ for the uncolored vertices}\label{fig-exIlist}
\end{subfigure}
\begin{subfigure}{0.2\textwidth}
\centering
\epsfig{file=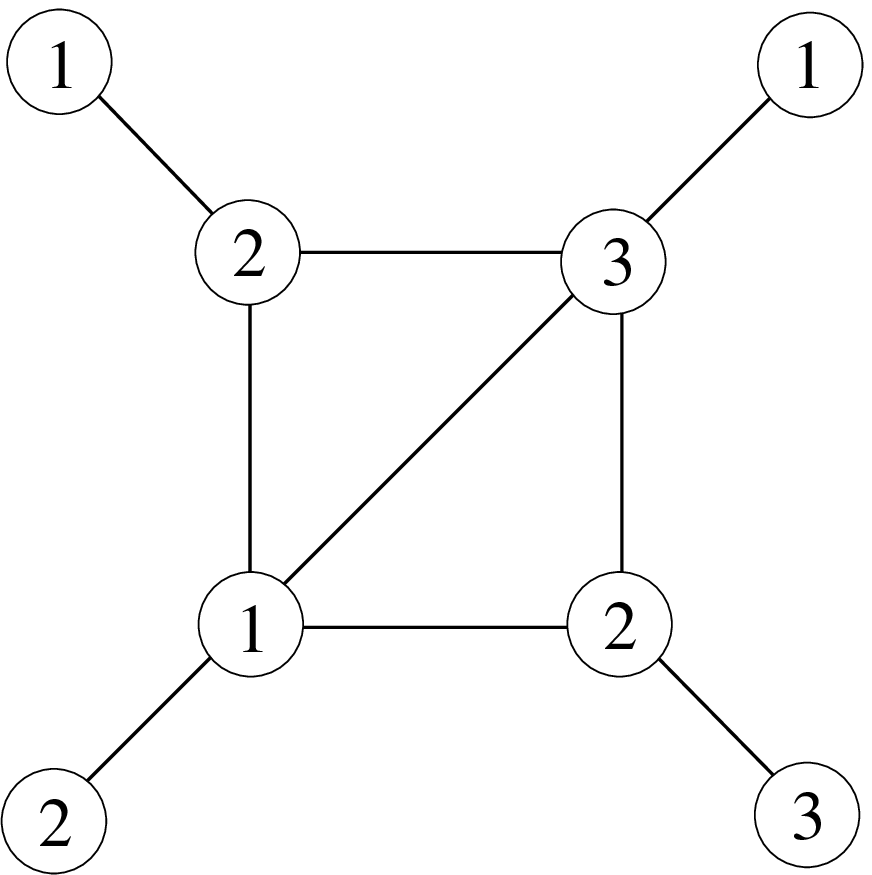, width=3.6cm}
\caption{Coloring $\mathscr C$}\label{fig-exII}
\end{subfigure}
\caption{A Sudoku coloring and its unique extension.}\label{f4}
\end{figure}
\nt The argument is similar to Example~\ref{ex-ex1}. We omit here. \rsq
\end{example}

\nt We now present a few results on list coloring of paths and cycles, which are used in the next section.

\begin{lemma}\label{lem-listcoloring-path} Let $L(x_i)$ be a list of colors of a vertex $x_i$ in the path   $P_n=x_1x_2\cdots x_n$, $1\le i\le n$. If
$|L(x_i)|\ge 2$ for each $i$, then there are at least two list colorings of $P_n$.
\end{lemma}

\begin{proof}
We color $x_i$ in natural order starting at $i=1$. When $x_i$, $1\le i\le n-1$, has been colored, then there is at least one color in $L(x_{i+1})$ can be assigned to $x_{i+1}$. So we have a list coloring for $P_n$. Since there are two choices for coloring $x_1$, we have the lemma.
\end{proof}



\begin{lemma}\label{lem-listcoloring-cycle} Suppose that there exists a list coloring of the cycle $C_n=x_1x_2\cdots x_nx_1$ such that the list of colors for each vertex $x_i$ satisfies the following conditions:
\begin{enumerate}[1.]
\item  $|L(x_i)|\ge 2$;
\item $L(x_i)\subseteq \{1,2,3\}$.
\end{enumerate}
Then there are at least two list colorings of $C_n$.
\end{lemma}

\begin{proof} To show this lemma, we may assume a weaker condition that $|L(x_i)|= 2$ for each $i$. Suppose $L(x_i)$ are the same for all $i,$ $1\le i\le n$. Since $C_n$ admits a list coloring, it follows that $n$ is even.  Now we are going to show that there exist two list colorings.

\nt Now we assume that not all $L(x_i)$ are same. Then there exist two adjacent vertices whose lists are different. Hence we may assume without loss of generality that $L(x_1)=\{1,2\}$ and $L(x_n)=\{1,3\}.$ Let $k$ be the least positive integer such that $2\leq k\leq n$ and $L(x_k)\neq L(x_1).$ If $k=n,$ then we assign color 3 to $x_1$ and by Lemma~\ref{lem-listcoloring-path} we get two list colorings for the path $C_n-x_n.$ Suppose $k\leq n-1.$ Let $\phi$ be a list coloring of $C_n.$ Let $P=x_1x_2\cdots x_{k-1}$ and $Q=x_kx_{k+1}\cdots x_n.$ We claim that there exists a list coloring of $C_n$ different from $\phi.$  We consider two cases.

\begin{enumerate}[A:]
\item Let $L(x_k)=\{1,3\}.$
\begin{enumerate}[{A}1.]
\item Suppose $\phi(x_n)=\phi(x_k)=3.$ Then by interchanging the colors of the vertices of $P,$ we get another list coloring of $C_n.$
\item Suppose $\phi(x_n)=\phi(x_k)=1.$ Then $P$ is an odd path and $\phi(x_1)=\phi(x_{k-1})=2.$ Now recolor $x_k$ with color 3. By Lemma~\ref{lem-listcoloring-path} we get a list coloring of the path $Q$ and this gives a list coloring of $C_n$ different from $\phi.$
\item Suppose $\phi(x_n)$ and $\phi(x_k)$ are distinct. Without loss of generality, let $\phi(x_n)=1$ and $\phi(x_k)=3.$ Then $\phi(x_1)=2$ and $\phi(x_{k-1})=\left\{\begin{array}{ll} 2 &\quad\mbox{if $k-1$ is odd}\\ 1 &\quad\mbox{if $k-1$ is even.}\end{array}\right.$ \\
If $k-1$ is odd we recolor $x_k$ by color 1 and by Case~A2, we get a list coloring of $C_n$ different from $\phi.$ If $k-1$ is even, we swap the colors 1 and 2 for $P$ and then recolor $x_k$ by 3. Similarly we recolor the inverse path $Q^{-1}=x_nx_{n-1}\cdots x_k$ and obtain a list coloring different from $\phi.$
\end{enumerate}

\item $L(x_k)=\{2,3\}.$

If $\phi(x_k)=\phi(x_n)=3,$ the proof is similar to that of Case A1. Suppose $\phi(x_k)\neq \phi(x_n).$

\begin{enumerate}[B1.]
\item Suppose $\phi(k_n)=1.$ Using Lemma~\ref{lem-listcoloring-path}, we get a list coloring of the path $R=x_nx_{n-1}\cdots x_1$ in which $x_n$ is assigned color 3. Since the color assigned for $x_1$ is 1 or 2, this coloring is also a list coloring of $C_n$ different from $\phi.$
\item Suppose $\phi(x_n)=3$ and $\phi(x_k)=2.$ Consider a list coloring of the path $x_kx_{k+1}\cdots x_nx_1x_2\cdots x_{k-1}$ in which $x_k$ is colored 3. Since the color assigned to $x_{k-1}$ is 1 or 2, this is also a list coloring of $C_n$ which is different from $\phi.$
\end{enumerate}
\end{enumerate}
Thus there exist two list  colorings of $C_n.$
\end{proof}

\nt We now proceed to present basic properties of Sudoku coloring.

\begin{observation} Let $G$ be a connected graph with $\chi(G)\geq 3$. Let $\mathscr C_0=\{V_1,V_2,\dots,V_r\}$ be a Sudoku coloring defined on an induced subgraph $G[S]$, for some $S\subset V(G)$. If $r\leq k-2,$ then at least two  new color classes are created in the extension process and hence the extension is not unique. Therefore $r=k-1$ or $k$.
\end{observation}

\nt Let $\mathscr C_0$ be an extendable coloring of $G[S]$ and let $W=V(G)\setminus S.$ Let $w_1\in W$ be the first vertex to be colored in the extension process. Let $\mathscr C_1$ be the coloring on $G[S\cup \{w_1\}].$ Let $w_2$ be the next vertex to be colored and let $\mathscr C_2$ be the resulting coloring of $G[S\cup \{w_1,w_2\}].$ Continuing this process we get a sequence of colorings $\mathscr C_0$, $\mathscr C_1$, $\mathscr C_2$, $\dots$, $\mathscr C_t$ where $t=|W|$ and $\mathscr C_t$ is the required extension of $\mathscr C_0$ to a $k$-coloring of $G.$ In the extension process, we update the color-list of the vertices as follows. If a vertex $w$ is assigned color $i,$ then $L(w)$ is replaced by the empty set $\varnothing$ and the color $i$ is removed from the list of all neighbors of $w.$

\begin{defi} Let $G$ be a graph with $\chi(G)=k\geq 3.$ Let $\mathscr C_0=\{V_1,V_2,\dots,V_{k-1}\}$ be an extendable coloring of $G[S]$ and let $w\in W=V\setminus S.$ The vertex $w$ is called a {\it color dominating vertex} if $w$ is adjacent to at least one vertex in each color class $V_i, 1\leq i\leq k-1.$
\end{defi}

\begin{observation} If $w$ is a color dominating vertex, then in the extension process, $k$ is the unique color assigned to $w.$
\end{observation}

\begin{defi} Let $G$ be a graph with $\chi(G)=k\geq 3.$ Let $\mathscr C_0=\{V_1,V_2,\dots,V_k\}$ be an extendable coloring of $G[S]$ and let $w\in V\setminus S.$ The vertex $w$ is called a {\it near-color dominating vertex} if there exists a unique color class $V_j$ such that $v$ is not adjacent to any vertex in $V_j.$
\end{defi}

\begin{observation} If $w$ is a near-color dominating vertex, then in the extension process, $j$ is the unique color that can be assigned to $w.$
\end{observation}

\begin{example} Consider the Sudoku coloring $\mathscr C_0$ of the graph in Fig.~\ref{f3}. There is no color or near-color dominating vertex. \rsq
\end{example}

\begin{example}  Let $C_{13}=v_1v_2\cdots v_{13}v_1$ be the $13$-cycle. Let $\mathscr C_0=\{V_1,V_2,V_3\}$ be a coloring of $C_{13}[S]$, where $S=\{v_i\;|\; 1\le i\le 13, i \mbox{ is odd}\}$,  $V_1=\{v_1,v_5,v_9\}, V_2=\{v_3,v_7,v_{11}\}$ and $V_3=\{v_{13}\}$. Here every vertex $w$ in $V \setminus S$ is a near-color dominating vertex and can be assigned a unique color. Thus $\mathscr C_0$ is a Sudoku coloring of $C_{13}$ and the unique extension $\mathscr C$ is given by $\mathscr C=(\{V_1\cup \{v_{12}\}, V_2,V_3\cup \{v_2,v_4,v_6,v_8,v_{10}\})$. \rsq

\end{example}

\begin{example} Consider the partial coloring $\mathscr C_0$ for the graph $C_6(K_4)$ given in Fig.~\ref{f7}.
\begin{figure}[H]
\centering
\epsfig{file=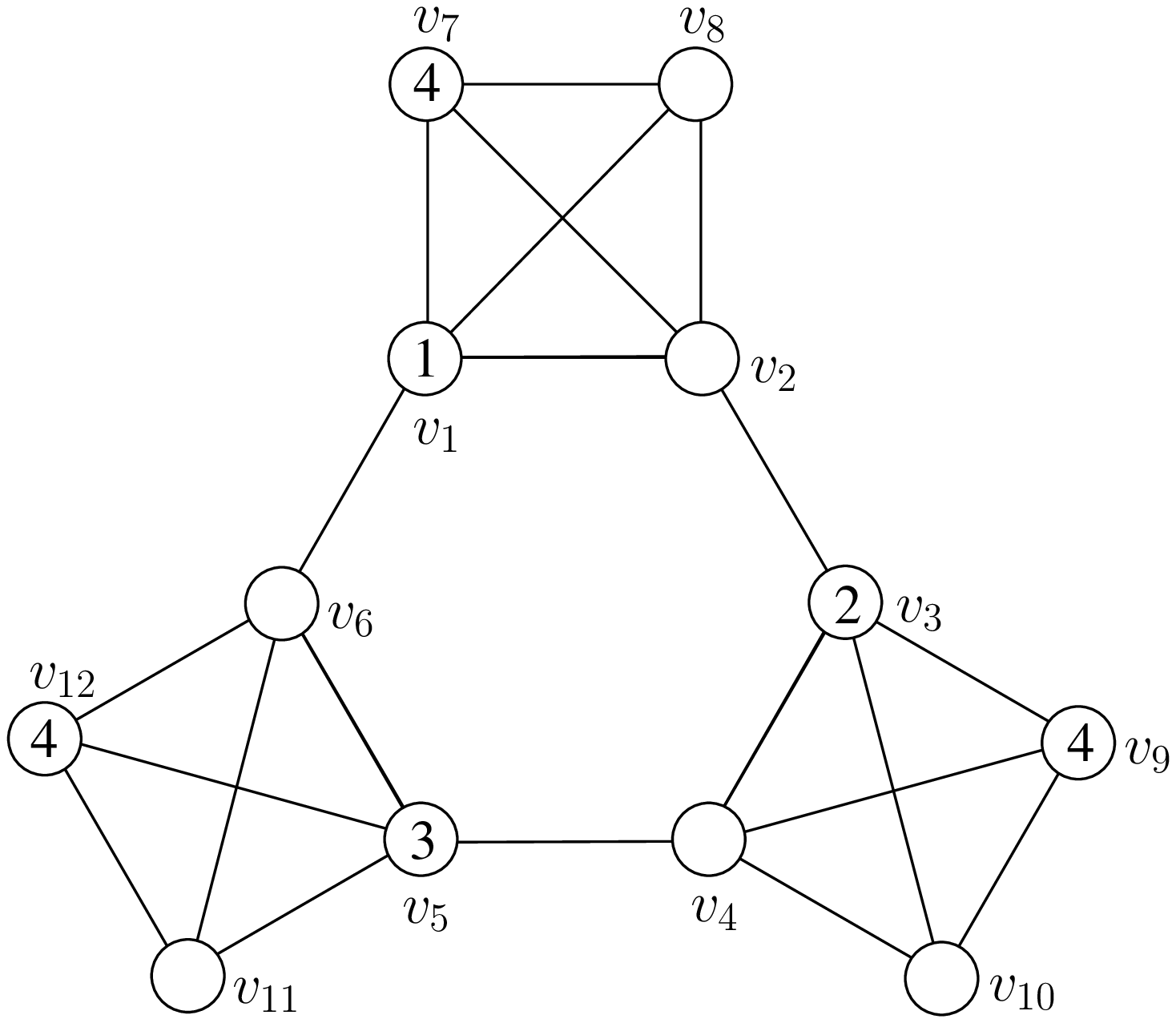, width=5cm}
\caption{A Sudoku coloring for $C_6(K_4)$.}\label{f7}
\end{figure}

\nt Clearly $\mathscr C_0=\{\{v_1\},\{v_3\},\{v_5\},\{v_7,v_9,v_{12}\}\}.$ The vertex $v_6$ is a near-color dominating vertex for $\mathscr C_0$ and receives the color 2. Hence $\mathscr C_1=\{\{v_1\},\{v_3,v_6\},\{v_5\},\{v_7,v_9,v_{12}\}\}.$ Now $v_{11}$ is a near-color dominating vertex for $\mathscr C_1.$ (Note that $v_{11}$ is not a near-color dominating vertex for $\mathscr C_0$). The unique color for $v_{11}$ is 1. Hence $\mathscr C_2=\{\{v_1,v_{11}\},\{v_3,v_6\}, \{v_5\},\{v_7,v_9,v_{12}\}\}.$ Proceeding like this we get
\begin{align*}
\mathscr C_3&=\{\{v_1,v_{11}\},\{v_3,v_6\},\{v_5,v_2\},\{v_7,v_9,v_{12}\}\},\\
\mathscr C_4&=\{\{v_1,v_{11}\},\{v_3,v_6,v_8\}, \{v_5,v_2\}, \{v_7,v_9,v_{12}\}\},\\
\mathscr C_5&=\{\{v_1,v_{11},v_4\},\{v_1,v_6,v_8\},\{v_5,v_2\},\{v_7,v_9,v_{12}\}\}\\
{\rm and}\ \mathscr C_6&=\{\{v_1,v_{11},v_4\},\{v_1,v_6,v_8\},\{v_5,v_2,v_{10}\},\{v_7,v_9,v_{12}\}\}
\end{align*}
$\mathscr C_6$ is the unique extension of $\mathscr C_0$ to a 4-coloring of $G.$ Also $W=V\setminus S=\{v_6,v_{11},v_2,v_8,v_4,v_{10}\}$ where the vertices in $W$ are listed in the order in which colors are assigned in the extension process. \rsq
\end{example}

\nt The above examples illustrate the significant role of color dominating vertex and near-color dominating vertex in the extension process of a Sudoku coloring. The following concept provides another tool in this regard.

\begin{defi} Let $G$ be a graph with $\chi(G)=k\geq 3.$ Let $\mathscr C_0=\{V_1,V_2,\dots,V_k\}$ be an extendable coloring defined on $G[S]$  for some $S\subset V$. Suppose $L$ is the corresponding list of the vertex set $V\setminus S$.  The vertex $w\in V\setminus S$ is called an {\it $i$-attractive vertex} if $w$ satisfies the following conditions.
\begin{enumerate}[(i)]
\item The chromatic number of the subgraph induced by the closed neighborhood $N[w]$ is $k.$
\item There exists a color $i\in \{1,2,\dots,k\}$ such that $i\notin L(v)$ for all $v\in N(w)$,  the neighborhood of $w$.
\end{enumerate}

\nt Clearly if $w$ is an $i$-attractive vertex, then in any extension of $\mathscr C_0$, the vertex $w$ is assigned color $i$. In solving a Sudoku puzzle, this condition is very often used in determining the entry in a cell. But there will be no such vertex in some cases. Please see Example~\ref{ex-exI}.
\end{defi}

\begin{lemma}\label{lem-pendant} Let $G$ be a graph with $\chi(G)=k\ge 3$. Suppose $\mathscr C_0$ is an extendable coloring of $G[S]$ for $S\subset V(G)$. If there is a pendant vertex $v\notin S$, then $\mathscr C_0$ is not a Sudoku coloring.
\end{lemma}
\begin{proof} Let $\mathscr C$ be an extension of $\mathscr C_0$ to $G-v$. Let $u$ be the unique neighbor of $v$ and let $i$ be the color assigned to $u.$ Then $v$ can be assigned any of the colors from the set $\{1,2,\dots,k\}\setminus \{i\}.$ Thus $\mathscr C_0$ has more than one extension and hence the result follows.
\end{proof}

\begin{lemma} \label{lem-K2} Let $G$ be a graph with $\chi(G)=k\ge 3$. Suppose $\mathscr C_0$ is an extendable coloring of $G[S]$ for $S\subset V(G)$. If there is an edge $xy$ for which $x,y\notin S$ such that $\deg(x)\le k-1$ and $\deg(y)\le k-1$, then $\mathscr C_0$ is not a Sudoku coloring of $G$.
\end{lemma}
\begin{proof} Let $\mathscr C$ be the extension of $\mathscr C_0$ to $G-\{x,y\}$. Now the lists of available colors of $x$ and $y$ are $\{i\;|\;1\le i\le k\}\setminus \{\mathscr C(u)\;|\; u\in N(x)\setminus\{y\}\}$ and $\{i\;|\;1\le i\le k\}\setminus \{\mathscr C(v)\;|\; v\in N(y)\setminus\{x\}\}$, respectively. By assumption, both of these lists contain at least $2$ colors. Thus, $\mathscr C_0$ can be extended to at least 2 different colorings for $G$.
\end{proof}

\section{Sudoku Numbers for Some Graphs}

\begin{theorem}\label{thm-bipart} Let $G$ be a connected graph of order $n.$ Then $sn(G)=1$ if and only if $G$ is a bipartite graph.
\end{theorem}

\begin{proof} Let $G$ be a connected bipartite graph. Since the bipartition $(V_1,V_2)$ is uniquely determined, if one vertex of $G$ is assigned color 1, then the 2-coloring of $G$ is uniquely determined. Thus, $sn(G)=1$.

\ms\nt Conversely, suppose $sn(G)=1.$ Then there exists a vertex $v$ of $G$ such that the 1-coloring $\mathscr C$ of $K_1=\{v\}$ determines a $\chi$-coloring of $G.$ If $\chi(G)\geq 3,$ then $\chi(G-v)\geq 2$ and any extension of $\mathscr C$ to a $\chi$-coloring of $G$ is not unique, which is a contradiction. Thus $\chi(G)=2$ and $G$ is bipartite.
\end{proof}

\begin{corollary} Every tree $T$ has $sn(T)=1$. \end{corollary}

\nt In what follows, we only consider graph $G$ with $\chi(G)\ge 3$.

\begin{theorem}\label{thm-tpart} Let $G$ be a uniquely colorable graph and let $\chi(G)=t.$ Then $sn(G) = t-1$. \end{theorem}

\begin{proof} Let $\mathscr C=\{V_1,V_2,\dots,V_t\}$ be the unique $t$-coloring of $G.$ Let $u_i\in V_i$ and let $S=\{u_1,u_2,\dots,u_{t-1}\}.$ Let $\mathscr C_0$ be a coloring  of $G[S]$ in which $u_i$ is assigned the color $i.$ Clearly $\mathscr C_0$  is uniquely extendable to a $t$-coloring of $G.$ Hence $sn(G)\leq t-1.$

\ms\nt Now, let $S\subseteq V(G)$ and $|S|\leq t-2.$ Clearly $\chi(G-S)\geq 2.$ Hence for any coloring of $G-S,$ its extension to a $\chi$-coloring of $G$ is not unique. Hence $sn(G)\geq t-1$ and the theorem follows.
 \end{proof}

\begin{corollary}\label{cor-Kn} Let $G$ be a complete $t$-partite graph with $t\geq3.$ Then $sn(G) = t-1$. \end{corollary}

\begin{theorem} For $n\ge 3$, $$sn(C_n)= \begin{cases} 1 & \mbox{ if } $n$ \mbox{ is even}, \\ \frac{n+1}{2} & \mbox{ otherwise.}\end{cases}$$ \end{theorem}

\begin{proof} The case $n$ is even follows from Theorem~\ref{thm-bipart}. So we assume that $n$ is odd. Hence $\chi(C_n) = 3$.

\nt The case $n=3$ is obvious. For $n\ge 5$, let $C_n=v_1v_2\cdots v_nv_1$. Let $S=\{v_j\;|\; 1\le j\le n-2, \mbox{ $j$ is odd}\}\cup\{v_{n-1}\}$.
Define a coloring $\mathscr C$ for $C_n[S]$ by  $\mathscr C(v_i)=1$ for $i\equiv 1\pmod{4}$, $\mathscr C(v_j)=2$ for $j\equiv3\pmod{4}$ and $\mathscr C(v_{n-1})=3$. This can be uniquely extended to a 3-coloring of $C_n$ since we must have $\mathscr C(v_k)=3$ for even $k\not= n-1$ and $\mathscr C(v_n)=2$. Thus, $\mathscr C$ is a Sudoku coloring of $C_n$. Since $C_n[S]$ is of order $\frac{n+1}{2}$, we have $sn(C_n)\le \frac{n+1}{2}$.

\ms\nt Suppose there is an extendable coloring $\phi$ of $G[S]$ with $|S| \le \frac{n-1}{2}.$ Then $|V(G)\setminus S|\ge \frac{n+1}{2}$. By pigeonhole principle, there is a $K_2$ subgraph which satisfies the condition of Lemma~\ref{lem-K2}. Thus $\phi$ is not a Sudoku coloring. Hence $sn(C_n)=\frac{n+1}{2}$.
\end{proof}

\nt For $m\ge 2$ and $1\le i\le m$, let $G_i$ be a simple graph with an induced subgraph $H$. An {\it amalgamation} of $G_1,\ldots, G_m$ over $H$ is the simple graph obtained by identifying the vertices of $H$ of each $G_i$ so that the obtained new graph contains a subgraph $H$ induced by the identified vertices. This subgraph $H$ is called the {\it common core} of the amalgamation of $G_1,\ldots, G_m$. Suppose $G$ is a graph with a proper subgraph $K_r$, $r\ge 1$. Let $A(mG, K_r)$ be the amalgamation of $m\ge 2$ copies of $G$ over $K_r$. Note that there may be many non-isomorphic $A(mG, K_r)$ graphs. When $r=1$, the graph is also known as one-point union of graphs. Note that $A(mK_2,K_1)\cong K_{1,m}$ and $A(mK_3, K_1)$ is the {\it friendship graph} $f_m$, $m\ge 2$.

\begin{theorem}\label{thm-Amal} For $m\ge 2$, $n\ge 3$ and $1\le r<n$, if $G=A(mK_n,K_r)$, then $sn(G) =m(n-r-1)+r-1$. \end{theorem}

\begin{proof} Note that $|V(G)| = r + m(n-r)$ and $\chi(G)=n$. Let $K\cong K_r$ be the common core of $G$ and let $G_i$ be the $i$-th copy of $K_n$.

\ms\nt Choose one vertex in $K$, say $x_0$, and choose one vertex in each $G_i-K$, say $x_i$. Let $S=V(G)\setminus \{x_i\;|\; 0\le i\le m\}$. Then $|S|=m(n-r-1)+r-1$.

\ms\nt Now color the vertices of $K-x_0$ from $1$ to $r-1$; those of $G_1-x_1$ from $r+1$ to $n-1$ and those of $G_i-x_i$ ($2\le i\le m$) from $r+2$ to $n$, arbitrary. Note that, if $n=r+1$, then we do not perform the last two assignments. Clearly, this coloring can be extended to an $n$-coloring of $G$ uniquely. Thus, $sn(G) \le m(n-r-1)+r-1$.

\ms\nt Suppose there is an extendable coloring $\phi$ of $G[S]$ with $|S| < m(n-r-1)+r-1$. Then $|V(G)\setminus S|\ge m+2$. Considering the $m+1$ subgraphs, $K$ and $G_i-K$, $1\le i\le m$, by pigeonhole principle, there is an edge $xy$ in either $K$ or $G_i-K$ for $i$, for which $x$ and $y$ have not been colored.

\ms\nt Let $\hat\phi$ be an extension of $\phi$. Then $\hat\phi(x)\ne \hat\phi(y)$. For convenience, we assume $\hat\phi(x)=1$ and $\hat\phi(y)=2$.

\begin{enumerate}[1.]
\item Suppose $x,y\in V(G_i-K)$ for some $i$. Since $\chi(G_i)=n$, $\{\hat\phi(v)\;|\; v\in V(G_i)\setminus\{x,y\}\}=\{j\;|\; 3\le j\le n\}$. Now, if we swap the colors of $x$ and $y$, then we get another extension of $\phi$ for all $i$. Hence $\phi$ is not a Sudoku coloring.
\item Suppose $x,y\in V(K)$. Similar to Case 1, we will obtain that $\{\hat\phi(v)\;|\; v\in V(G_i)\setminus\{x,y\}\}=\{j\;|\; 3\le j\le n\}$ for all $i$. Again, we may swap the colors of $x$ and $y$ to obtain another extension of $\phi$. Hence $\phi$ is not a Sudoku coloring.
\end{enumerate}
Thus, $sn(G) = m(n-r-1)+r-1$.
\end{proof}

\begin{corollary}  For $m\ge 2$, $sn(f_m)=m$.  \end{corollary}

\ms\nt A {\it tadpole} graph $T(n,m)$ is obtained from a cycle $C_n$ and a $P_m$ by identifying a vertex of $C_n$ to an end vertex of $P_m$.

\begin{theorem} For $n\ge 3$ and $m\ge 2$, $$sn(T(n,m))=\begin{cases} 1 &\mbox{ if } n \mbox{ even,} \\ \lfloor\frac{n+m}{2}\rfloor & \mbox{ if } n \mbox{ odd}. \end{cases}$$ \end{theorem}

\begin{proof} Let $C_n=v_1v_2\cdots v_nv_1$ and $P_m =u_1u_2\cdots u_m$ with $v_1 = u_1$. Also let $G=T(n,m)$. The case $n$ is even follows from Theorem~\ref{thm-bipart}. So, following we assume $n$ is odd.

\ms \nt Suppose $m$ is even. Let $S=\{v_i\;|\; 2\le i\le n-1, \mbox{$i$ is even}\}\cup\{u_j\;|\; 2\le j\le m, \mbox{$j$ is even} \}$ Define a coloring $\mathscr C$ on $G[S]$ by
\[\mathscr C(u_j)=\begin{cases}2 &  \mbox{ if }j\equiv2\pmod{4}\\ 3 & \mbox{ if } j\equiv 0\pmod{4}\end{cases} \mbox{ and } \mathscr C(v_i)=\begin{cases}3 &  \mbox{ if } i\equiv2\pmod{4}\\ 2 & \mbox{ if } i\equiv 0\pmod{4}\end{cases}\]

\ms\nt Now, we can extend this coloring $\mathscr C$ to a unique 3-coloring of $T(n,m)$ with $\mathscr C(v_n)=3,2$ depending on whether $n\equiv 1$ or $3\pmod 4$, and $\mathscr C(v)=1$ for each remaining vertex $v$. Thus, $sn(T(n,m))\le \frac{n+m-1}{2}=\lfloor\frac{n+m}{2}\rfloor$.

\ms\nt Suppose there is an extendable coloring $\phi$ of $G[S]$ with $|S| \le \frac{n+m-3}{2}.$ Then $|V(G)\setminus S|\ge \frac{n+m+1}{2}$. By Lemma~\ref{lem-pendant} we may assume that $v_m\in S$. Now, consider $\frac{n+m-1}{2}$ subsets $\{u_{2j}, u_{2j+1}\}$, $1\le j\le (m-2)/2$; $\{v_{2i}, v_{2i+1}\}$, $1\le i\le (n-1)/2$; and $\{v_1\}$. By pigeonhole principle, there is a $K_2$ in $G- S$ which satisfies the condition of Lemma~\ref{lem-K2}. Thus $\phi$ is not a Sudoku coloring. Thus $sn(T(n,m))= \lfloor\frac{n+m}{2}\rfloor$.

\ms\nt Suppose $m$ is odd. Let $S=\{v_i\;|\; 2\le i\le n-1, \mbox{$i$ is even}\}\cup\{u_j\;|\; 1\le j\le m, \mbox{$j$ is odd} \}.$ Define a coloring $\mathscr C$ on $G[S]$ by
\[\mathscr C(u_j)=\begin{cases}1  &  \mbox{ if }j\equiv 1\pmod{4}\\ 3 & \mbox{ if } j\equiv 3\pmod{4}\end{cases} \mbox{ and } \mathscr C(v_i)=\begin{cases}3 &  \mbox{ if } i\equiv2\pmod{4}\\ 2 & \mbox{ if } i\equiv 0\pmod{4}\end{cases}\]

\ms\nt Now, we can extend this coloring $\mathscr C$ to a unique 3-coloring of $T(n,m)$ with $\mathscr C(v_n)=3,2$ depending on whether $n\equiv 1$ or $3\pmod 4$, , $\mathscr C(v_i)=1$ for each remaining vertex $v_i$ and $\mathscr C(u_j)=2$ for each remaining vertex $u_j$. Thus, $sn(T(n,m))\le \frac{n+m}{2}=\lfloor\frac{n+m}{2}\rfloor$.

\ms\nt Suppose there is an extendable coloring $\phi$ of $G[S]$ with $|S| \le \frac{n+m-2}{2}$, then $|V(G)\setminus S|\ge \frac{n+m}{2}$. By Lemma~\ref{lem-pendant} we may assume that $v_m\in S$. Now, consider $\frac{n+m}{2}-1$ edges $u_{2j-1}u_{2j}$, $1\le j\le (m-1)/2$; and $v_{2i}v_{2i+1}$, $1\le i\le (n-1)/2$. By pigeonhole principle, there is an edge $xy$ in $G- S$. If $xy\ne u_1u_2$, then by Lemma~\ref{lem-K2} $\phi$ is not a Sudoku coloring. Now let $xy=u_1u_2$. We may additionally assume that  no edge of the graph $G-u_1$ lies in $G-S$.

\ms\nt  If $v_n\notin S$, then the path $v_nu_1u_2$ lies in $G-S$. We extend $\phi$ to $G-\{v_n, u_1, u_2\}$. Then the color-list of $v_n, u_1, u_2$ are of at least 2 colors. By Lemma~\ref{lem-listcoloring-path}, $\phi$ can be extended to at least two $3$-colorings of $G$. Thus $\phi$ is not a Sudoku coloring.

\ms\nt If $v_n\in S$, then $v_2\notin S$ by the additional condition. This is the same case as the above case.

\ms\nt Thus, $sn(T(n,m))= \frac{n+m}{2}$.
\end{proof}

\ms\nt A {\it lollipop} graph $L(n,m)$ is obtained from a $K_n$ and a path $P_m$ by identifying a vertex of $K_n$ to an end vertex of $P_m$. Note that $T(3,m)=L(3,m)$.

\begin{theorem}\label{thm-lollipop} For $n\ge 4$ and $m\ge 2$, $sn(L(n,m)) = n+m-3$. \end{theorem}

\begin{proof} Let the vertices of $K_n$ be $v_i, 1\le i\le n$ and $P_m =u_1u_2\cdots u_m$ with $v_1=u_1$. Obviously $\chi(L(n,m)) = n$.

\ms\nt Let $G=L(n,m)$. Let $S=\{v_i\;|\; 3\le i\le n\}\cup\{u_j\;|\; 2\le j\le m\}$. Then $|S|=n+m-3$. Define $\mathscr C(v_i)=i$, $3\le i\le n$, $\mathscr C(u_j)=2$ for even $j$ and $\mathscr C(u_j)=1$ for odd $j$, where $2\le j\le m$. If we extend $\mathscr C$, then $v_1=u_1$ must be colored by $1$ and $v_2$ by $2$. So $\mathscr C$ is a Sudoku coloring. Hence $sn(G)\le n+m-3$.

\ms\nt Let $\phi$ be an extendable coloring of $G[S]$, for some $S\subset V(G)$ with $|S|<n+m-3$. In other words $|V(G)\setminus S|\ge 3$. Let $x,y,z$  be three vertices in $V(G)\setminus S$. We extend $\phi$ to $G-\{x,y,z\}$ first.

\ms\nt Suppose $x=u_i$, $2\le i\le m$. Then the color-list of $x$ contains at least $n-2$ colors. Thus $\phi$ has at least two extensions.

\ms\nt Suppose $x,y,z$ are in $K_n$. Without loss of generality let $x=v_2$ and $y=v_3$. Hence $\phi$ is not a Sudoku coloring by Lemma~\ref{lem-K2}.

\ms\nt Thus $sn(G)= n+m-3$.
\end{proof}

\ms \nt Let $C_{2n}(K_m)$ be a graph obtained from a $2n$-cycle $C_{2n}$ by identifying every alternate edge of $C_{2n}$ with an edge of a distinct complete graph $K_m$. Thus, $C_{2n}(K_m)$ contains $n$ copies of complete subgraph $K_m$ such that the edges not belong to any $K_m$ are alternate edges of $C_{2n}$. If  $C_{2n}=x_1x_2\cdots x_{2n-1}x_{2n}x_1$, then $x_{2i-1}x_{2i}$ is an edge of the $i$-th complete graph $K_m$, $1\le i\le n$. We will denote the $i$-th complete graph $K_m$ by $K^i$. When we remove the vertices $x_{2i-1}$ and $x_{2i}$ from the $K^i$, then the resulting subgraph $H^i\cong K_{m-2}$.

\begin{theorem}  For $n\ge 2$, $m\ge 3$, $sn(C_{2n}(K_m)) = n(m-2)$. \end{theorem}

\begin{proof} Keep the notation defined above. Color any $m-3$ vertices of each $H^i$ by $\{j\;|\; 4\le j\le m\}$, $1\le i\le n$. For even $n$, color $x_{4j+1}$ by 1 and $x_{4j+3}$ by 2 for $0\le j\le (n-2)/2$. For odd $n$, color $x_{4j+1}$ by 1; $x_{4j+3}$ by 2 for $0\le j\le (n-3)/2$ and $x_{2n-1}$ by 3.

\ms\nt Clearly, this is a proper coloring for the subgraph of $C_{2n}(K_m)$ induced by the colored vertices. To extend this coloring to the whole graph, all $x_{2j}$'s must be colored by 1, 2, or 3 uniquely. After that, the color of the last uncolored vertex in each $H^i$ is also fixed.
Thus, $sn(C_{2n}(K_m))\le n(m-2)$.

\ms\nt Suppose there is an extendable coloring $\phi$ of $C_{2n}(K_m)[S]$ with $|S|\le n(m-2)-1$. Then at least $2n+1$ vertices of $C_{2n}(K_m)$ have not been colored under $\phi$. By pigeonhole principle, at least one $K^i$, say $i=1$ for convenience, contains at least three uncolored vertices, say $u, v, w$. We extend $\phi$ to be an $m$-coloring of $C_{2n}(K_m)-\{u,v,w\}$. Now, the color-lists of $u, v, w$ are the same and are of size 3. By Lemma~\ref{lem-listcoloring-cycle}, $\phi$ is not a Sudoku coloring. Thus $sn(C_{2n}(K_m))= n(m-2)$.\end{proof}

\begin{example} Two Sudoku colorings for $C_6(K_4)$ and $C_8(K_4)$ are given in Fig.~\ref{fn}.
\begin{figure}[H]
\centerline{\epsfig{file=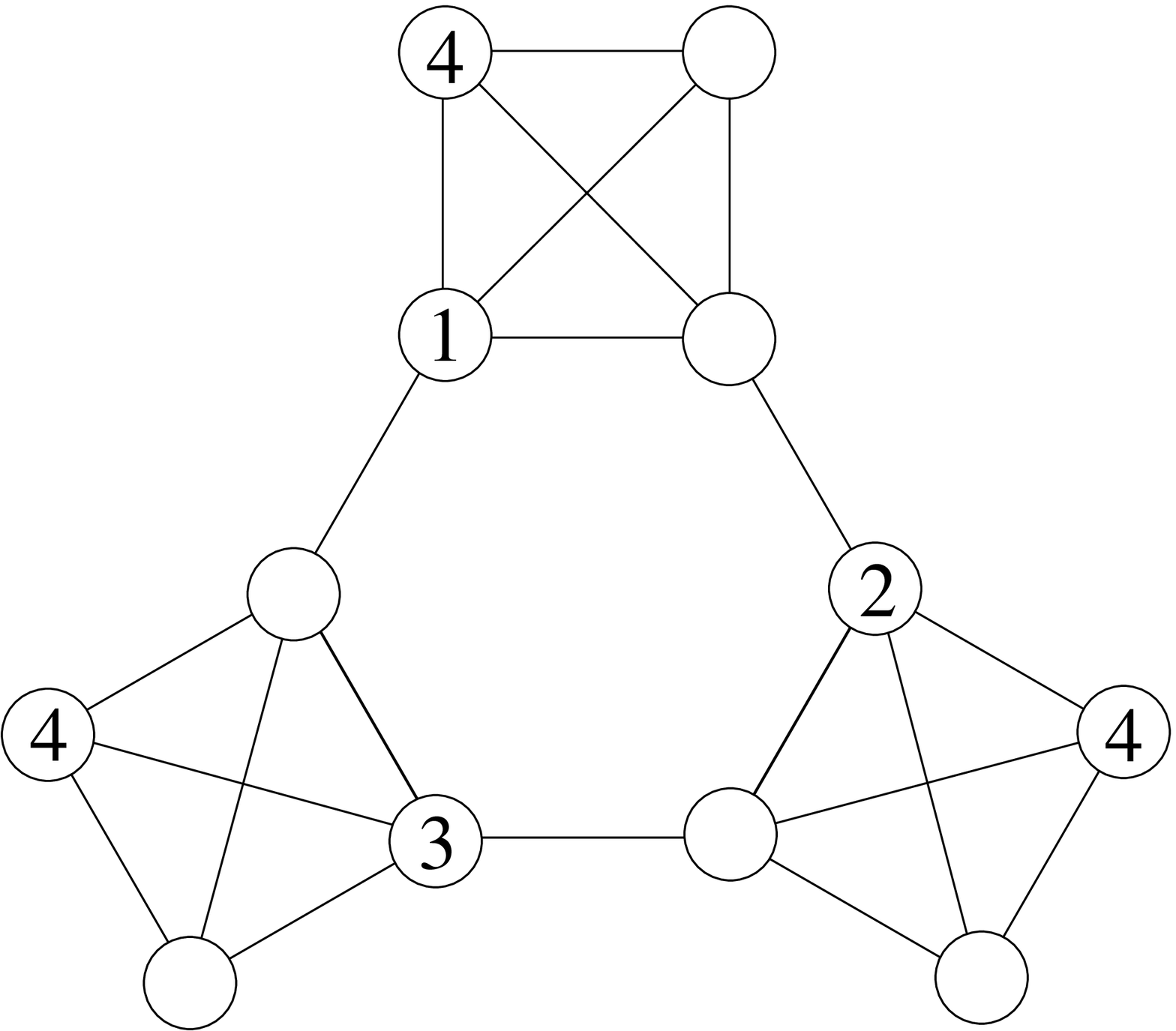, height=3cm}\qquad\qquad\epsfig{file=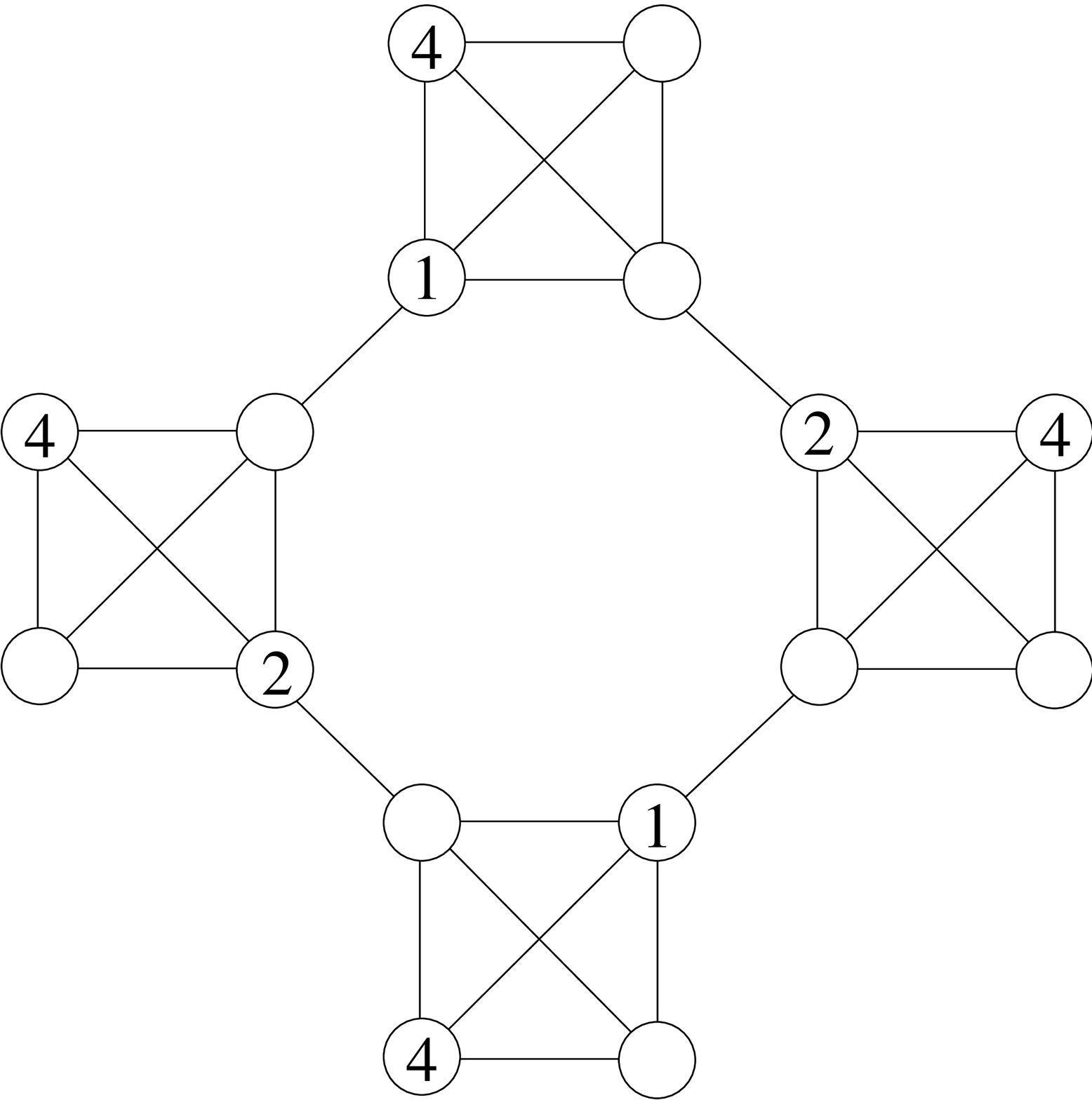, height=3cm}}
\caption{Sudoku colorings of $C_6(K_4)$ and $C_8(K_4)$}\label{fn}
\end{figure}
\rsq
\end{example}

\nt Keep the notation of defining the graph $C_{2n}(K_m)$. Let $C_{2n}(K_m^-)$ be the graph obtained from $C_{2n}(K_m)$ by removing
the $n$ edges $x_{2i-1}x_{2i}$, $1\le i\le n$. Thus $C_{2n}(K_m^-)$ is an $(m-1)$-regular graph.

\begin{theorem}  For $n\ge 2$ and $m\ge 4$, $sn(C_{2n}(K_m^-)) = (m-3)n+1$. \end{theorem}

\begin{proof} Let $G=C_{2n}(K_m^-)$. Let $y_{i,j}$, $1\le j\le m-2$, be the vertices of $H^i$, $1\le i\le n$. Note that $\chi(G)=m-1$. We shall first show that $sn(G)\le (m-3)n+1$.

\ms\nt Suppose $n$ is even. Define a partial coloring $\mathscr C$ on $C_{2n}(K_m^-)$ by $\mathscr C(y_{i,j})=j+1$ for $1\le i\le n$ and $2\le j\le m-2$, $\mathscr C(y_{1,1})=2$. To extend $\mathscr C$ as an $(m-1)$-coloring, the odd path $x_2x_3y_{2,1} x_4  \cdots  x_{2n-1} y_{n,1} x_{2n} x_1$ must be colored alternatively by 1,2 starting from $x_2$ and ending at $x_1$. So $\mathscr C$ is a Sudoku coloring for $G$.

\ms \nt Suppose $n$ is odd. Define a partial coloring $\mathscr C$ on $C_{2n}(K_m^-)$ by $\mathscr C(y_{i,j})=j+1$ for $1\le i\le n-1$ and $2\le j\le m-2$, $\mathscr C(y_{1,1})=2$, $\mathscr C(y_{n,j})=j+1$ for $3\le j\le m-2$ (no this case if $m=4$), and $\mathscr C(y_{n,2})=1$. To extend $\mathscr C$ as an $(m-1)$-coloring, the even path $x_2 x_3 y_{2,1} x_4 \cdots  x_{2n-1} y_{n,1} x_{2n} x_1)$ must be colored  alternatively by 1,2 starting from $x_2$ and ending at $x_{2n-3}$. And then the last 6 vertices $y_{n-1,1}, x_{2n-2}, x_{2n-1}, y_{n,1}, x_{2n}, x_1$ must be colored by 1, 2, 3, 2, 3,1,  respectively. So $\mathscr C$ is a Sudoku coloring for $G$.

\ms\nt So we have $sn(C_{2n}(K_m^-))\le (m-3)n+1$ for both cases of $n$.

\ms\nt Suppose there is an extendable coloring $\phi$ of $G[S]$ with $|S|\le (m-3)n$. Let the extension of $\phi$ be $\psi$. Then at least $3n$ vertices of $G$ have not been colored under $\phi$.
\begin{enumerate}[1.]
\item Suppose $y_{i,j_1}$ and $y_{i,j_2}$ have not been colored under $\phi$, for some $i$, $1\le j_1< j_2\le m-2$. Then $\psi(y_{i,j_1})\ne \psi(y_{i,j_2})$. Now we may exchange the colors of $y_{i,j_1}$ and $y_{i,j_2}$ to obtain another $(m-1)$-coloring for $G$. Thus $\phi$ is not a Sudoku coloring for $G$.

\item Suppose at most one $y_{i,j}$ is uncolored by $\phi$, for each $i$. In other word, at most 3 vertices of each $K^i-x_{2i-1}x_{2i}$ have not been colored. By pigeonhole principle, exactly three vertices of $K^i-x_{2i-1}x_{2i}$ have not been colored for each $i$. That is, $x_{2i-1}$, $y_{i,j_i}$ and $x_{2i}$ have not been colored, for some $j_i$. All uncolored vertices induce a $3n$-cycle $C$. The color-list of each vertex containing exactly two colors. By Lemma~\ref{lem-listcoloring-cycle}, $\phi$ is not a Sudoku coloring for $G$.
\end{enumerate}
Thus, $sn(C_{2n}(K_m^-)) = (m-3)n+1$.
\end{proof}

\begin{example}A Sudoku coloring for $C_{10}(K_5^-)$ and its extension are given in Fig.~\ref{fig-C10K5-}.
\begin{figure}[H]
\centerline{\begin{tabular}{m{4.5cm}m{1.5cm}m{4.5cm}}
\epsfig{file=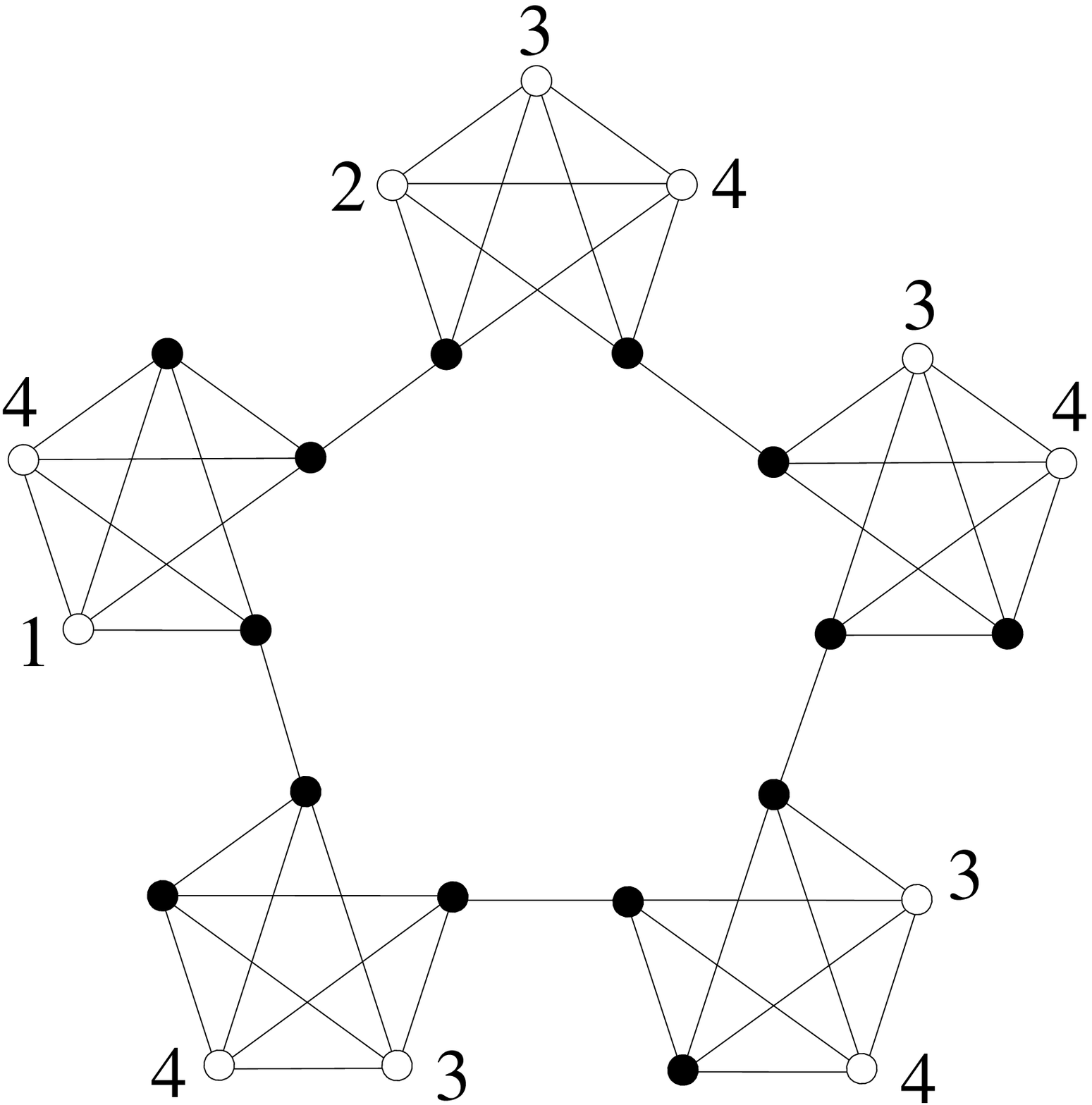, width=4.5cm} & \quad\Large $\rightarrow$ &\epsfig{file=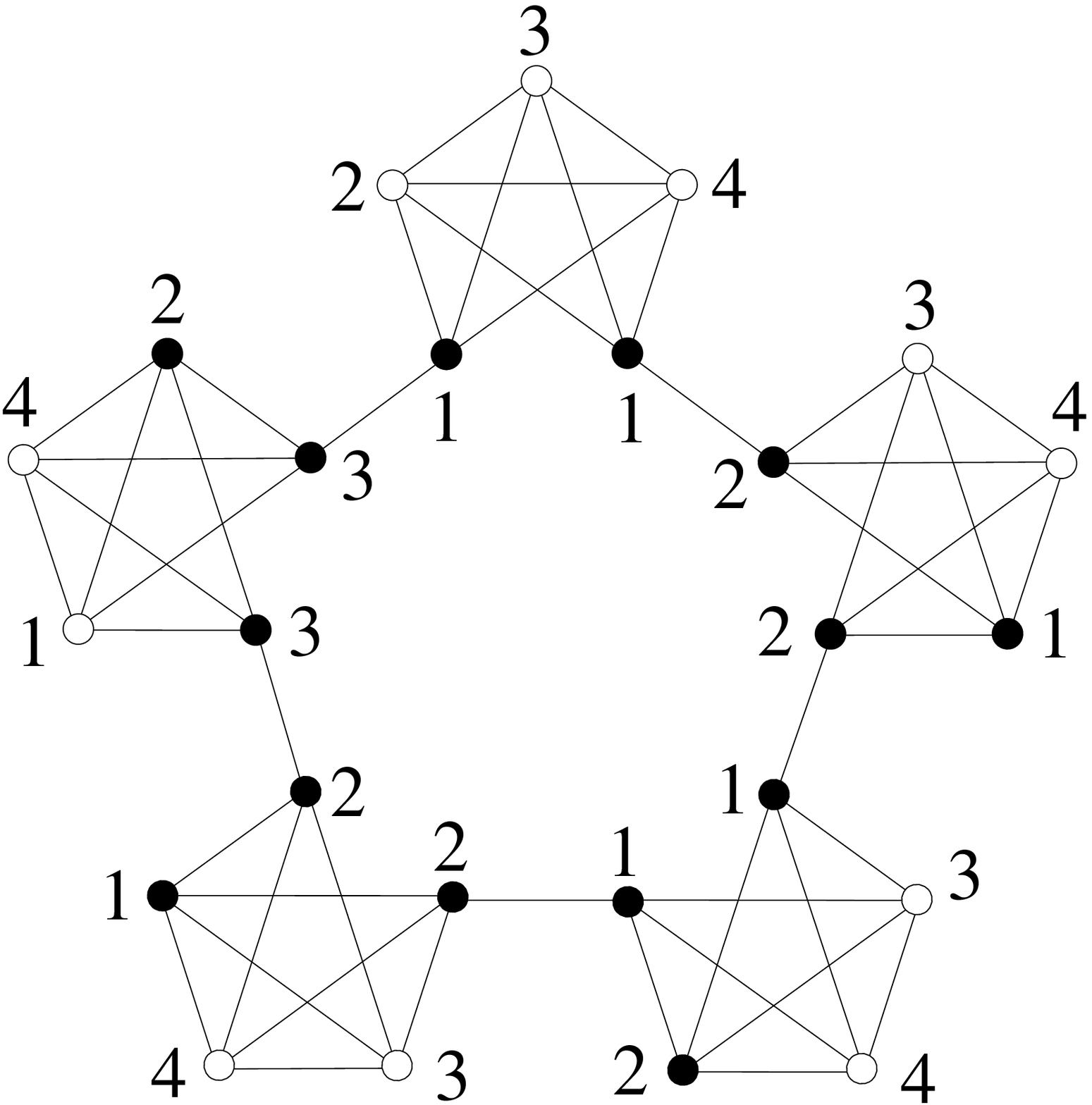, width=4.5cm}\end{tabular}}
\caption{A Sudoku coloring and its extension for $C_{10}(K_5^-)$}\label{fig-C10K5-}
\end{figure}
\rsq
\end{example}

\nt Let $G_1 = K_3$ with vertices $x_1,y,z$. For $i\ge 2$, let $G_i$ be obtained from $G_{i-1}$ by adding a vertex $x_i$ adjacent to two adjacent vertices of $G_{i-1}$ (to form a new triangle) so that $G_i$ is a $2$-connected plane graph whose faces are triangles. Note that for $i\ge 3$, $G_i$ is not unique. Consequently, every maximal outerplanar graph of order $n+2$ is a possible $G_n$, $n\ge 1$. A {\it fan graph} $F_n$, $n\ge 2$, is obtained from $P_n=v_1v_2\cdots v_n$ by joining a vertex $u$ to every vertex of $P_n$. Thus, $F_n$ is a possible $G_{n-1}$.

\begin{theorem} For $i\ge 1$, let $G_i$ be defined as above. Then $sn(G_i)=2$.   \end{theorem}

\begin{proof} Since $\chi(G_i)=3$, we have $sn(G_i)\ge 2$ for $i\ge 1$. Let $S=\{y,z\}$ and $G_i[S]$ has a coloring $\mathscr C(y)=1, \mathscr C(z)=2$. This coloring can be extended to a unique $3$-coloring of $G_i$ since we must color $x_1$ by $3$ and for $i\ge 2$, we must color $x_i$ by exactly one of $1$, $2$, or $3$ which is different to the colors assigned to its two neighbors.  So, $sn(G_i)\le 2$. Thus, $sn(G_i)=2$.     \end{proof}

\begin{corollary} Every fan graph (and maximal outerplanar graph) $G$ has $sn(G)=2$. \end{corollary}

\nt  A {\it wheel graph} $W_n$, $n\ge 3$ of order $n+1$ is obtained from an $n$-cycle $C_n=v_1v_2\cdots v_nv_1$ by joining a vertex $u$ to every vertex of $C_n$. In each of these graphs, the maximum degree vertex is called the {\it core} of the graph.

\begin{theorem} For $n\ge 3$, $sn(W_n)=2$ if $n$ is even, $sn(W_3)=3$ and $sn(W_n) = \frac{n+1}{2}$ if $n\ge 5$ is odd. \end{theorem}

\begin{proof} Since $W_n$ is not bipartite, $sc(W_n)\ge 2$.

\ms\nt Suppose $n$ even. It is known that $\chi(W_n)=3$. Let $S=\{u,v_1\}$, then $W_n[S]$ has a coloring $\mathscr C(u)=1$, $\mathscr C(v_1) = 2$. This coloring can be uniquely extended to a 3-coloring of $W_n$ since we must have $\mathscr C(v_i) = 2$ for remaining odd $i$, and $\mathscr C(v_i)=3$ for even $i$. So, $sn(W_n)\le 2$. Hence $sn(W_n)= 2$.

\ms\nt Suppose $n$ is odd. It is known that $\chi(W_n)=4$. The case $W_3 = K_4$ follows from Corollary~\ref{cor-Kn}. We assume $n\ge 5$. Let $S = \{v_i\,|\,1\le i\le n, \mbox{ $n$ is odd}\}$. Now, $W_n[S]$ has a coloring $\mathscr C(v_n)=4$, $\mathscr C(v_i)=2$ for $i\equiv 1\pmod{4}, i\ne n$ and $\mathscr C(v_i)=3$ for $i\equiv 3\pmod{4}, i\ne n$. This coloring can be uniquely extended to a 4-coloring of $W_n$ since we must have $\mathscr C(u)=1$, $\mathscr C(v_i)=4$ for even $i \le n-3$, and $\mathscr C(v_{n-1})=2$ if $n\equiv 1\pmod{4}$ while $\mathscr C(v_{n-1}) = 3$ if $n\equiv 3\pmod{4}$. So, $sn(W_n)\le \frac{n+1}{2}$.

\ms\nt Suppose there is an extendable coloring $\phi$ of $W_n[S]$ with $|S| \le\frac{n-1}{2}$. Since $C_n$ is Hamiltonian, the number of components of $C_n-(S\setminus\{u\})$ is most $|S|\le \frac{n-1}{2}$. Since $C_n-(S\setminus\{u\})$ contains $\frac{n+1}{2}$ vertices, there is a $K_2$ subgraph in $C_n-(S\setminus\{u\})$. By Lemma~\ref{lem-K2}, $\phi$ is not a Sudoku coloring. Hence, $sn(W_n) = \frac{n+1}{2}$.
\end{proof}

\section{Conclusion and Scope}

For any connected graph $G$ of order $n,$ we have $1\leq sn(G)\leq n-1.$ We have shown that $sn(G)=1$ if and only if $G$ is bipartite. Further $sn(K_n)=n-1.$ In this context we propose the following conjecture.

\begin{conjecture} For a graph $G$ of order $n,$ $sn(G)=n-1$ if and only if $G=K_n$. \end{conjecture}

\nt From the algorithmic perspective we propose the following conjecture.

\begin{conjecture} Given a connected graph $G$ with $X(G)=k\geq 3$ and a partial coloring $\mathscr C_0$ defined on $G[S]$ where $S\subset V,$ the problem of deciding whether $\mathscr C_0$ is an extendable coloring is NP-complete.

\ms\nt  If $G$ is disconnected graph with $c$ components, then $sn(G)\geq c.$
\end{conjecture}

\begin{problem} Study the Sudoku number of disconnected graphs.  \end{problem}


\begin{thebibliography}{99}

\bibitem{Bondy} J.A. Bondy, U.S.R. Murty, {\it Graph theory with applications}, New York, MacMillan, 1976.
%

\bibitem{Erdos} P. Erd\H{o}s, A.L. Rubin, and H. Taylor, Choosability in graphs, {\it  Congr. Numer.}, {\bf  26}, 125--157 (1979).

\bibitem{Vizing} V.G. Vizing, Vertex colorings with given colors, {\it Metody Diskret. Analiz.} (in Russian), {\bf 29}, 3--10, (1976).

\bibitem{Sudoku} https://abcnews.go.com/blogs/headlines/2012/06/can-you-solve-the-hardest-ever-sudoku,\\ accessed Sept. 24, 2021.
\end{thebibliography}
\end{document}